\documentclass[11pt,oneside,a4paper]{article}
\usepackage{amsmath}
\numberwithin{equation}{section}
\usepackage{amsfonts}
\usepackage{tikz}
\usepackage{tikz-cd}
\usetikzlibrary{positioning, calc}
\usepackage{amssymb}
\usepackage{amsthm}
\usepackage{dsfont}
\usepackage{enumitem}
\usepackage{hyperref}
\usepackage{indentfirst}
\usepackage{mathrsfs}
\usepackage[top=1in,bottom=1in,left=1.25in,right=1.25in]{geometry}

\date{}
\allowdisplaybreaks

\begin{document}

\renewcommand{\baselinestretch}{1.2}
\renewcommand{\arraystretch}{1.0}

\title{\bf Ore extensions of multiplier Hopf coquasigroups}
\author
{\textbf{Rui Zhang, Na Zhang, Yapeng Zeng}
\footnote{College of Sciences, Nanjing Agricultural University, Nanjing 210095, China.\\
 E-mail: 2024111007@stu.njau.edu.cn, 23123110@stu.njau.edu.cn,
  23123127@stu.njau.edu.cn},
    \textbf{Tao Yang} \footnote{Corresponding author.
        E-mail: tao.yang@njau.edu.cn}
}
\maketitle

\begin{center}
    \begin{minipage}{12.cm}
        \textbf{Abstract}
        In this paper, Ore extensions of multiplier Hopf coquasigroups are studied. Necessary and sufficient conditions for the Ore extension of a regular multiplier Hopf coquasigroup to be a multiplier Hopf coquasigroup are given. Furthermore, the isomorphism between two such Ore extensions is discussed.
        \\

        \\

    \end{minipage}
\end{center}
\normalsize

\section{Introduction}
\def\theequation{\thesection.\arabic{equation}}
\setcounter{equation}{0}

 Hopf algebras have evolved into increasingly general formalisms over the past decades. 
 In the early 1990s, Van Daele \cite{V94} introduced the concept of multiplier Hopf algebras, a framework that dispenses with the unitality requirement intrinsic to classical Hopf algebras and thereby furnishes new tools for the study of quantum groups and non-commutative algebras. This theory has since been systematically expanded to include an intrinsic duality theory, a comprehensive representation theory, and applications in quantum field theory.  Hopf coquasigroups have undergone substantial development, particularly through the foundational work in 2010 of Klim and Majid \cite{KM}, who formulated a weak form of coassociativity. Their approach establishes profound links with algebraic topology and opens new avenues for the study of non-(co)associative symmetries. 
 
 The Ore extension is a well-known method in ring extensions used to construct noncommutative rings and algebras. 
 Moreover, within the theory of Hopf algebras, Ore extensions are significant for constructing noncommutative and non-cocommutative Hopf algebras, as demonstrated in the works of Panov \cite{P03}, Beattie \cite{B00, BDG00}, and Nenciu \cite{N01, N}.
 Subsequently, researchers investigated the feasibility of endowing these generalized Hopf-type structures with Ore extensions (see Figure 1). Zhao and Lu \cite{ZL12} initiated the systematic study of Ore extensions in the context of multiplier Hopf algebras, while Jiao and Meng \cite{JM14} formulated Ore extensions for Hopf coquasigroups and established the notion of Hopf coquasigroup–Ore extensions.

\begin{figure}[htbp]
\centering
\begin{tikzpicture}[node distance=1cm and 0.25cm, auto]
  \node (start) [align=center]
  {{Ore Extensions of}\\ {Hopf algebras}};
  
  \node (middle1) [above right=of start, align=center] 
    {{Ore Extensions of}\\ {Hopf coquasigroups}};
  
  \node (middle2) [below right=of start, align=center] 
    {{Ore Extensions of}\\{multiplier Hopf algebras}};
  
  \node (end) [right=0.01mm of $(middle1.east)!0.5!(middle2.east)$, align=center] 
    {{Ore Extensions of}\\ {multiplier Hopf coquasigroups}};

  \draw[->, thick] (start) -- 
    node[above left=2mm, align=center] {\footnotesize Jiao \& Meng\\ \footnotesize \cite{P03}} (middle1);
  
  \draw[->, thick] (start) -- 
    node[below left=2mm, align=center] {\footnotesize Zhao \& Lu\\ \footnotesize \cite{ZL12} } (middle2);
  
  \draw[->, thick] (middle1) -- 
    node[above=2mm, sloped] {\footnotesize ?} (end);
  
  \draw[->, thick] (middle2) -- 
    node[below=2mm, sloped] {\footnotesize ?} (end);

\end{tikzpicture}
\caption{Diagram of the Ore Extension Structure Relationships} 
 \label{fig:ore_extension}
\end{figure}
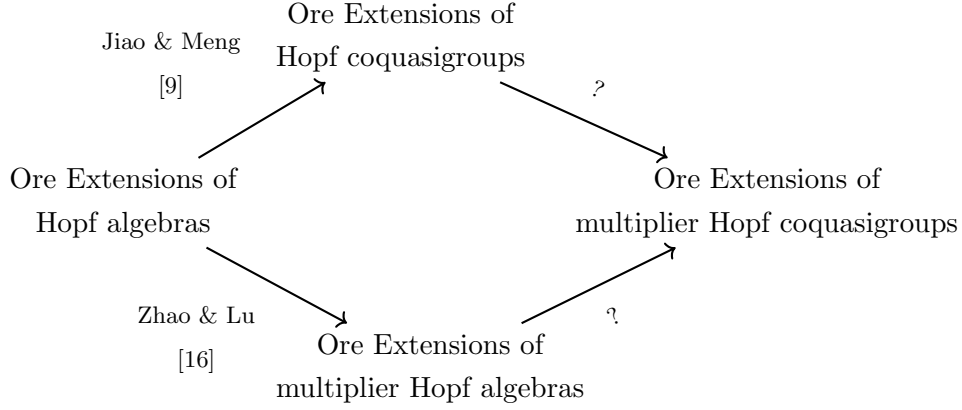

 In 2022, Yang \cite{Y22} established a unifying framework for multiplier Hopf algebras and Hopf coquasigroups through the introduction of multiplier Hopf coquasigroups. This theoretical advancement not only bridged these two important algebraic structures but also provided a solution to the biduality problem for a class of infinite-dimensional Hopf quasigroups. 
 Motivated by these recent advances, the present work addresses the following fundamental question: Can Ore extensions be  realized within the framework of multiplier Hopf coquasigroups?

 This paper provides a positive answer, and its structure is as follows. 
 In Section 2, we provide an introduction to key concepts: Ore extensions, multiplier algebras, and multiplier Hopf coquasigroups, which will be referenced in subsequent sections.

 Section 3 discusses Ore extensions within the framework of multiplier Hopf coquasigroups. Additionally, one of the main results of this study (Theorem 3.4) is presented, offering a definitive condition for Ore extensions of regular multiplier Hopf coquasigroups to also be regular. This theorem affirms the previous question. Meanwhile, we consider the $*$-algebra version.

 Moving on to Section 4, we discuss the isomorphism for Ore extensions and establish the criteria under which two such extensions are considered isomorphic.

\section{Preliminaries}
\def\theequation{\thesection.\arabic{equation}}
\setcounter{equation}{0}

 Throughout this paper, all linear spaces under consideration are defined over a fixed field $k$ (for example, the field of complex numbers $ \mathds {C} $).

\subsection{Multiplier Hopf coquasigroups}

 Let $A$ be an (associative) algebra over $k$, with or without an identity element, but with a nondegenerate product. This condition implies that if $a \in A$ and $ab = 0$ for all $b \in A$ or $ba = 0$ for all $b \in A$, then $a$ must be equal to $0$. 

 Recall from \cite{V94, V98, VV25} that a left multiplier of $A$ is a linear map $\rho_1: A \rightarrow A$ satisfying   $\rho_1(ab) = \rho_1(a)b$ for all $ a,b \in A $. Similarly, a right multiplier of $A$ is a linear map $\rho: A \rightarrow A$ such that $\rho_2(ab) = a\rho_2(b)$ for all $a, b \in A$. 
 A multiplier of $A$ consists of a pair $ (\rho_1, \rho_2) $ comprising a left and a right multiplier, where $\rho_2(a)b = a\rho_1(b)$ holds for all $a,b \in A$. 
 For any algebra $A$, there exist a left multiplier algebra, a right multiplier algebra, and a multiplier algebra of $A$, denoted as $L(A)$, $R(A)$, and $M(A)$, respectively. 
 
 Specifically, $M(A)$ is characterized as the largest algebra with identity containing $A$ as an essential two-sided ideal.
 Moreover, in the case where $x \in M(A) $ and $xa = 0$ for all $a \in A$ or $ax = 0$ for all $a \in A$, it follows that $x = 0$. Additionally, we examine the tensor algebra $A \otimes A$, which remains nondegenerate, and we have its multiplier algebra $M(A \otimes A)$.
 There are natural embeddings
 $$ A \otimes A \subseteq M(A) \otimes M(A) \subseteq M(A \otimes A). $$
 Generally, when $A$ lacks an identity element, these two embeddings are strict. In the scenario where $A$ possesses an identity element $1_{A}$, the product is evidently nondegenerate, and we observe that $M(A) = A$ and $M(A \otimes A) = A \otimes A$.

 Let $A$ and $B$ be nondegenerate algebras. If a homomorphism $f: A \longrightarrow M(B)$ is nondegenerate (meaning $f(A)B = B$ and $ Bf(A) = B $), then it can be uniquely extended to a homomorphism $M(A) \longrightarrow M(B)$, also denoted by $f$.

 A \emph{multiplier Hopf coquasigroup}, as defined in \cite{Y22}, is a nondegenerate associative algebra $A$ equipped with algebra homomorphisms $\Delta: A \longrightarrow M(A \otimes A)$ (comultiplication), $\varepsilon: A \longrightarrow k$ (counit), and a linear map $S: A \longrightarrow A$ (antipode). The key conditions for this structure are as follows:
\begin{enumerate}[nosep, label=(\arabic*)]
    \item  For any $a, b \in A$, both $T_{1}(a \otimes b) = \Delta(a)(1 \otimes b)$ and $T_{2}(a \otimes b) = (a \otimes 1) \Delta(b)$ belong to $A \otimes A$.
    \item The counit satisfies $(\varepsilon \otimes \iota)T_{1}(a \otimes b) = ab = (\iota \otimes \varepsilon)T_{2}(a \otimes b)$, where $\iota$ denotes the identity map, and we will use this notation henceforth.
    \item The antipode $S$ is defined to be antimultiplicative and anticomultiplicative, such that for any $a, b \in A$ the following hold:
          \begin{eqnarray}
              & (m \otimes \iota)(S \otimes \Delta)\Delta(a) = 1_{M(A)} \otimes a = (m \otimes \iota)(\iota \otimes S \otimes \iota)(\iota \otimes \Delta)\Delta(a), \label{2.1}\\
              & (\iota \otimes m)(\Delta \otimes S)\Delta(a) = a \otimes 1_{M(A)}  = (\iota\otimes m)(\iota \otimes S \otimes \iota)(\Delta \otimes \iota)\Delta(a). \label{2.2}
          \end{eqnarray}
\end{enumerate}
In the context of a multiplier Hopf coquasigroup $(A, \Delta)$, it is termed \emph{regular} if the antipode $S$ is bijective. The maps $T_1$ and $T_2$ are often called the \emph{Galois maps}.

   Multiplier Hopf coquasigroups, as introduced in \cite{Y22}, extend the concept of Hopf coquasigroups to a non-unital scenario. In contrast to the Hopf coquasigroups discussed in \cite{K, KM}, the algebraic structure of a multiplier Hopf coquasigroup $A$ may not necessarily possess an identity element in general. Moreover, the comultiplication map $\Delta$ does not merely map elements from $A$ to $A \otimes A$, but rather it maps them into the multiplier algebra $M(A \otimes A)$. This extension weakens the coassociativity property typically found in multiplier Hopf algebras while also harmonizing the concepts of multiplier Hopf algebras and Hopf coquasigroups.

 To denote the comultiplication, we employ the adapted Sweedler notation as introduced in \cite{V08}: $\Delta(a) = a_{(1)} \otimes a_{(2)}$ for $a \in A$.

\subsection{Ore extension}

\subsubsection*{Ore extension of an algebra}

 We start by recalling the definition of the Ore extension of an algebra, as introduced in \cite{K95}. Given an algebra $A$ with an identity element $1_A$ and an algebra endomorphism $\tau$ of $A$, a linear endomorphism $\delta$ of $A$ is referred to as a $\tau$-derivation if for all $a, b \in A$,
\begin{align*}
\delta(ab) = \delta(a)b + \tau(a) \delta(b).
\end{align*}
Moreover, it is a consequence of the above condition that $\delta(1_{A}) = 0$.

The Ore extension $A[y; \tau, \delta]$ of an algebra $A$ is constructed as an algebra that extends $A$ with the defining relation
\begin{align*}
ya = \tau(a)y + \delta(a)
\end{align*}
for all $a \in A$.

\subsubsection*{Ore extension of Hopf algebras}

 Let $A$ and $A[y; \tau, \delta]$ denote Hopf algebras over $k$. The Hopf algebra $A[y; \tau, \delta]$ is termed a Hopf-Ore extension if there exists an element $r \in A$ such that $\Delta(y) = y \otimes 1 + r \otimes y$, and $A$ is a Hopf subalgebra of $A[y; \tau, \delta]$.

In accordance with Theorem 1.3 in \cite{P03}, the conditions necessary and sufficient for the Hopf algebra $A[y; \tau, \delta]$ to qualify as a Hopf-Ore extension are as follows:

$(1)$ There exists a character $\chi: A \rightarrow k$ satisfying 
\begin{align*}
    \tau(a) = \chi(a_{(1)})a_{(2)}  \qquad \text{for all $ a \in A $};
\end{align*}

$(2)$ The following condition holds
\begin{align*}
    \chi(a_{(1)})a_{(2)}=Ad_r(a_{(1)})\chi(a_{(2)});
\end{align*}

$(3)$ The $\tau$-derivation $\delta$ fulfills the relation
\begin{align*}
    \Delta\big(\delta(a)\big) = \delta(a_{(1)}) \otimes a_{(2)} + ra_{(1)} \otimes \delta(a_{(2)}).
\end{align*}

\subsubsection*{Ore extension of Multiplier Hopf algebras}

 By generalizing the Ore extension of Hopf algebras to the non-unital case, Zhao and Lu \cite{ZL12} introduced the Ore extension of multiplier Hopf algebras. 
 
 For more details, let $(A, \Delta)$ be a regular multiplier Hopf algebra, and consider $A[y; \tau, \delta]$, an Ore extension of the algebra $A$. When $\tau$ is a nondegenerate algebra endomorphism, the $\tau$-derivation $\delta$ can be uniquely extended to the multiplier algebra $ M(A)$, denoted as $\overline{\delta}: M(A) \rightarrow M(A)$.
 In this case, $$ A[y; \tau, \delta] \subseteq M(A)[y;\tau,\delta] \subseteq M(A [y;\tau,\delta] ). $$
  Where $\tau$ is surjective, according to Proposition 2.3 of Section 2 in \cite{ZL12}, the product of $A[y; \tau, \delta]$ is nondegenerate. 
  
  Zhao and Lu first constructed $M(A)[y; \tau, \delta]$, an Ore extension of $M(A)$. The comultiplication, counit, and antipode of $A$ can be extended to $M(A)[y; \tau, \delta]$:
  $\overline{\Delta}: M(A)[y; \tau, \delta] \rightarrow M(A[y; \tau, \delta] \otimes A[y; \tau, \delta])$,  $\overline{\varepsilon}: M(A)[y; \tau, \delta] \rightarrow k$, $\overline{S}: M(A)[y; \tau, \delta] \rightarrow M(A[y; \tau, \delta])$.
 By restricting these extensions to the subspace $A[y; \tau, \delta]$, the authors obtained the comultiplication, counit, and antipode on the Ore extension $A[y; \tau, \delta]$ of $A$:
 \begin{align*}
 & \overline{\Delta} \big|_{A[y; \tau, \delta]}: A[y; \tau, \delta] \rightarrow M(A[y; \tau, \delta]\otimes A[y; \tau, \delta]),\\
 & \overline{\varepsilon}\big|_{A[y; \tau, \delta]}: A[y; \tau, \delta] \rightarrow k,\\
 & \overline{S}\big|_{A[y; \tau, \delta]}: A[y; \tau, \delta] \rightarrow M(A[y; \tau, \delta]).
 \end{align*}
 These maps were respectively denoted by $\Delta'$, $\varepsilon'$, and $S'$.

 The pair $(A[y; \tau, \delta], \Delta')$ is called a Hopf-Ore extension of the multiplier Hopf algebra $A$, if $(A[y; \tau, \delta], \Delta', \varepsilon', S')$ is a regular multiplier Hopf algebra such that $ \overline{\Delta}(y) = y \otimes 1_{M(A)} + r \otimes y $ for some $r \in M(A)$, where $1_{M(A)}$ is the identity element of $M(A)$.

 For a regular multiplier Hopf algebra $A$ with an Ore extension $A[y; \tau, \delta]$, where $\tau$ is surjective, $A[y; \tau, \delta]$ forms a Hopf-Ore extension if and only if there exists a group-like element $r \in M(A)$ satisfying the following equivalent conditions:

$(1)$ There exists a character $\chi:A \rightarrow k$ such that, for any $a \in A$, the following equality holds:
\begin{align*}
    \tau(a)=(\chi \otimes \iota)\Delta(a)=r\big((\iota \otimes \chi) \Delta(a)\big)r^{-1}.
\end{align*}

$(2)$ For the $\tau$-derivation $\delta$, for all $a \in A$, the relation
\begin{align*}
    \Delta \big(\delta(a)\big) = (\delta \otimes \iota)\Delta(a) + (r\otimes 1_{M(A)})\big((\iota \otimes \delta)\Delta(a)\big)
\end{align*}
is satisfied.

 As described in Section 2.1 of this paper, we will refer to these unique extensions as $\overline{\Delta}$, $\overline{S}$, $\overline{\varepsilon}$ simply as $\Delta$, $S$, $\varepsilon$ in the following without causing any confusion.

\subsubsection*{Ore extension of Hopf coquasigroups}

 Let $A$ and $A[y; \tau, \delta]$ denote Hopf coquasigroups. As defined in \cite{JM14}, the Hopf coquasigroup $A[y; \tau, \delta]$ is classified as a Hopf coquasigroups-Ore extension when $\Delta(y)=y \otimes 1 + r \otimes y$ for a certain $r \in A$, and $A$ is identified as a sub-Hopf coquasigroup of $A[y; \tau, \delta]$.

According to Theorem 3.3 in \cite{JM14}, a Hopf coquasigroup $A$ transforms into a Hopf coquasigroup-Ore extension denoted by $H = A[y, \tau, \delta ]$ if and only if the following conditions are met:

$(1)$ A character $\chi: A \rightarrow k$ exists such that for each $a \in A$, the condition $\tau(a) = \chi(a_{(1)})a_{(2)}$ is satisfied.

$(2)$The subsequent relations are valid:
\begin{align*}
    \chi(a_{(1)})a_{(2)(1)}\otimes a_{(2)(2)}
    = \text{Ad}_r(a_{(1)})\chi(a_{(2)(1)})\otimes a_{(2)(2)}
    = \chi(a_{(1)(1)})a_{(1)(2)}\otimes a_{(2)}.
\end{align*}

$(3)$ For the $\tau$-derivation $\delta$, the relation 
\begin{align*}
    \Delta(\delta(a)) = \delta(a_{(1)})\otimes a_{(2)}+ra_{(1)}\otimes\delta(a_{(2)})
\end{align*}
holds for all $a \in A$.

\section{Ore extensions of multiplier Hopf coquasigroups}
\def\theequation{\thesection.\arabic{equation}}
\setcounter{equation}{0}

 Let $(A, \Delta)$ denote a regular multiplier Hopf coquasigroup as defined in \cite{Y22}. The objective of this section is to establish a criterion for an Ore extension of a regular multiplier Hopf coquasigroup to also qualify as a regular multiplier Hopf coquasigroup.

 First, we shall introduce the definition of the Ore extension of regular multiplier Hopf coquasigroups, which extends the Ore extensions of Hopf coquasigroups and multiplier Hopf algebras as presented in \cite{JM14, ZL12}.
\\

\textbf{Definition 3.1.} Let $A$ and $A[y; \tau, \delta]$ be regular multiplier Hopf coquasigroups. The multiplier Hopf coquasigroup $(A[y; \tau, \delta], \Delta)$ is called an \emph{Ore extension of multiplier Hopf coquasigroup} (or an MHC-Ore extension in short) if the following equation holds:
\begin{equation}
     \Delta(y) = y \otimes 1_{M(A)} + r \otimes y  \label{3.1}
\end{equation}
for a group-like element $r \in M(A)$.
\\

 \textbf{Remark 3.2.}
 $(1)$ In what follows, we will simply write $1_{M(A)}$ as $1$ without causing any confusion. As stated in \cite{ZL12}, the element $y \otimes 1 + r \otimes y$ in equation $(\ref{3.1})$ belongs to $M(A[y; \tau, \delta] \otimes A[y; \tau, \delta])$. As usual, we denote $Ad_r(a) = raS(r) = rar^{-1}$.

 $(2)$ When the comultiplication $\Delta$ of the multiplier Hopf coquasigroup $A$ satisfies coassociativity, $A$ degenerates to a multiplier Hopf algebra and the Ore extension aligns with the one outlined in \cite{ZL12}. 
 As introduced in Section 2 of \cite{ZL12}, all extensions depend on the algebraic structure of $A$, then extensions on multiplier Hopf coquasigroups are similar.

 $(3)$ If the underlying algebra of the multiplier Hopf coquasigroup $A$
 possesses an identity element, $A$ can be classified as a Hopf coquasigroup, and the Ore extension corresponds to the one delineated in \cite{JM14}.
\\

Now, we present the properties associated with the Ore extension of a regular multiplier Hopf coquasigroup.
\\

\textbf{Proposition 3.3.} When $A[y; \tau, \delta]$ is an MHC-Ore extension of $A$, then
\begin{align*}
    &\varepsilon(y)=0, \quad S(y)=-r^{-1}y.
\end{align*}

\textbf{Proof.} By applying the counit map $\varepsilon \otimes \iota$ to both sides of equation $(\ref{3.1})$, we obtain 
\begin{align*}
    y&=(\varepsilon \otimes \iota)\Delta(y)\\
    &=(\varepsilon \otimes \iota)(y \otimes 1 + r \otimes y)\\
    &=\varepsilon(y)1 + \varepsilon(r)y\\
    &=\varepsilon(y)1 + y.
\end{align*}
Hence, it follows that $\varepsilon(y) = 0$.

Next, applying the same map $\varepsilon \otimes \iota$ again to equation $(\ref{3.1})$, we have
\begin{align*}
   \varepsilon(y) &=m(\iota \otimes S)\Delta(y)\\
                  &=m(\iota \otimes S)(y \otimes 1 + r \otimes y)\\
                  &=y + rS(y).
\end{align*}
Given that $\varepsilon(y) = 0$, we conclude that $S(y) = -r^{-1}y$. 
$\hfill \square$
\\

We now proceed to state the principal theorem of this section, providing necessary and sufficient conditions under which the Ore extension of a multiplier Hopf coquasigroup again admits a multiplier Hopf coquasigroup structure. 
\\

\textbf{Theorem 3.4.} Let $A$ and $A[y; \tau, \delta]$ be regular multiplier Hopf coquasigroups. Assume that $\tau$ is a surjective nondegenerate algebra endomorphism on $A$. Then the regular multiplier Hopf coquasigroup $A[y; \tau, \delta]$ is an MHC-Ore extension of $A$ if and only if there is a group-like element $r \in M(A)$ such that the following conditions hold:

$(1)$ There exists a character $\chi $: $ A \rightarrow k$ satisfying 
  \begin{equation}
    \tau(a)=\chi(a_{(1)})a_{(2)}=(\chi \otimes \iota)\Delta(a) \label{3.2}
  \end{equation}
for any $a \in A$.

$(2)$ The following equations hold
 \begin{equation}
    \chi (a_{(1)})a_{(2)(1)} \otimes a_{(2)(2)} = Ad_r(a_{(1)})\chi (a_{(2)(1)}) \otimes a_{(2)(2)} = \chi (a_{(1)(1)})a_{(1)(2)} \otimes a_{(2)} \label{3.3}
 \end{equation}
for any $a \in A$.

$(3)$The $\tau$-derivation $\delta$ satisfies the relation
\begin{equation}
    \Delta \big(\delta (a)\big) = (\delta \otimes \iota)\Delta(a) + (r\otimes 1)\big((\iota \otimes \delta)\Delta(a)\big)  \label{3.4}
\end{equation}
for any $a \in A$. 

\textbf{Proof.} The proof proceeds in two parts, we begin with the necessity.

\textbf{1. Proof of necessity.}

Analogous to the construction in \cite{ZL12}, the comultiplication $\Delta$ can be extended from $A$ to $A[y; \tau, \delta]$ by $\Delta(y) = y \otimes 1 + r \otimes y$. Then the algebra homomorphism $\Delta$ preserves the relation $ya = \tau(a)y + \delta(a)$ for any $a \in A$, i.e.,
\begin{align*}
\Delta(y)\Delta(a) = \Delta\big(\tau(a)\big)\Delta(y) + \Delta\big(\delta(a)\big).
\end{align*}

For any $b \in A$, we have
\begin{align*}
     &\Delta(y)\Delta(a)(1 \otimes b)\\ 
     &= (y \otimes 1 + r \otimes y)( a_{(1)}  \otimes a_{(2)}b) \\
     &=  ya_{(1)} \otimes a_{(2)}b + ra_{(1)} \otimes ya_{(2)}b \\
     &= \tau(a_{(1)})y \otimes a_{(2)}b + \delta(a_{(1)}) \otimes a_{(2)}b + ra_{(1)} \otimes \tau(a_{(2)}b)y + ra_{(1)} \otimes \delta(a_{(2)}b) \\
     &\overset{\bigstar}{=} \tau(a_{(1)})y \otimes a_{(2)}b + \delta(a_{(1)}) \otimes a_{(2)}b + ra_{(1)} \otimes \tau(a_{(2)})\tau(b)y+ r{a_{(1)}} \otimes \delta({a_{(2)}})b \\
     &\quad + r{a_{(1)}} \otimes \tau({a_{(2)}})\delta(b) \\
     &= \big(\tau(a_{(1)}) \otimes a_{(2)}\big)(1 \otimes b)(y \otimes 1) + \big(Ad_r(a_{(1)}) \otimes \tau(a_{(2)})\big)\big(r \otimes \tau(b)\big)(1 \otimes y) 
     \\
     &\quad + \big(Ad_r(a_{(1)}) \otimes \delta(a_{(2)})\big)(r \otimes b) + \big(Ad_r(a_{(1)}) \otimes \tau(a_{(2)})\big)\big(r \otimes \delta(b)\big) \\
     &\quad + \big(\delta(a_{(1)}) \otimes a_{(2)}\big)(1 \otimes b),
\end{align*}
where the identity $ra_{(1)} \otimes \tau(a_{(2)}b)y = ra_{(1)} \otimes \tau(a_{(2)})\tau(b)y$ at step $\bigstar$  follows from the fact that $\tau$ is a nondegenerate algebra homomorphism. Furthermore,
\begin{align*}
      &\Big(\Delta \big(\tau (a)\big)\Delta(y) + \Delta \big(\delta (a)\big)\Big)(1 \otimes b)\\
      &= \Delta \big(\tau (a)\big)(y \otimes b + r \otimes yb) + \Delta \big(\delta (a)\big)(1 \otimes b)\\
      &= \Delta \big(\tau (a)\big)(y \otimes b) + \Delta \big(\tau (a)\big)(r \otimes yb) + \Delta \big(\delta (a)\big)(1 \otimes b)\\
      &= \Delta \big(\tau (a)\big)(1 \otimes b)(y \otimes 1) + \Delta \big(\tau (a)\big)\big(r \otimes \tau (b)y\big) \\
      &\quad  + \Delta \big(\tau (a)\big)\big(r \otimes \delta (b)\big) + \Delta \big(\delta (a)\big)(1 \otimes b)\\
      &= \Delta \big(\tau (a)\big)(1 \otimes b)(y \otimes 1) + \Delta \big(\tau (a)\big)\big(r \otimes \tau (b)\big)(1 \otimes y)  \\
      &\quad + \Delta \big(\tau (a)\big)\big(r \otimes \delta (b)\big) + \Delta \big(\delta (a)\big)(1 \otimes b).
\end{align*}
Therefore, the equality $\Delta(y)\Delta(a) = \Delta(\tau(a))\Delta(y) + \Delta(\delta(a))$ holds if and only if the following three conditions are satisfied: for all $b\in A$,
\begin{align*}
    \Delta \big(\tau (a)\big)(1 \otimes b) &= \tau (a_{(1)}) \otimes a_{(2)}b,\\
    \Delta \big(\tau (a)\big)\big(r \otimes \tau (b)\big)&= ra_{(1)} \otimes \tau (a_{(2)}b),\\
    \Delta \big(\tau (a)\big)\big(r \otimes \delta (b)\big) + \Delta \big(\delta (a)\big)(1 \otimes b)&= ra_{(1)} \otimes \delta (a_{(2)}b) + \delta (a_{(1)}) \otimes a_{(2)}b.
\end{align*}
Using the properties of $\delta \otimes \iota$ and $\iota \otimes \delta$ from \cite[Lemma 2.5]{ZL12}, along with the extensions of $\tau \otimes \iota$, we reformulate the above conditions as follows:
\begin{eqnarray}
    &\Delta\big(\tau(a)\big) =(\tau \otimes \iota)\Delta(a)= \tau(a_{(1)}) \otimes a_{(2)}, \label{3.5}\\
    &\Delta\big(\tau(a)\big) = A{d_r}(a_{(1)}) \otimes \tau(a_{(2)}), \label{3.6}\\
    &\Delta \big(\delta (a)\big) = (\delta \otimes \iota)\Delta(a) + (r\otimes 1)\big((\iota \otimes \delta)\Delta(a)\big). \nonumber
\end{eqnarray}
The final condition is specified by equation $(\ref{3.4})$. 

We now demonstrate that equations $(\ref{3.5})$ and $(\ref{3.6})$ imply the existence of a character $\chi$ such that both $(\ref{3.2})$ and $(\ref{3.3})$ are valid. Following the approach in the proof of Proposition 3.4 in \cite{ZL12}, we define a map $\chi : A \rightarrow R(A)$ by 
\begin{align*}
    \chi (a)b = m( \iota \otimes S)\big(( b \otimes 1 )( \tau  \otimes \iota )\Delta (a) \big)
\end{align*}
for any $a,b \in A$. Since $\tau$ is nondegenerate, for any $b$ there exist elements $a_i$ ,$b_i$ such that $\sum_{i=1}^{n} {{a_i}} \tau \left( {{b_i}} \right) = b$. Hence, the expression above can be rewritten as:
\begin{align*}
\chi(a)b &= m( \iota  \otimes S )\Big( \sum_{i=1}^{n} \big( a_i \tau ( b_i ) \otimes 1 \big) \big( ( \tau \otimes \iota  ) \Delta(a) \big) \Big)\\
 &= m( \iota  \otimes S )\Big( \sum_{i=1}^{n} ( a_i \otimes 1 ) \big( \tau ( b_i ) \otimes 1 \big) \big( ( \tau \otimes \iota  ) \Delta(a) \big)\Big)\\
 &= m( \iota  \otimes S )\bigg( \sum_{i=1}^{n} ( a_i \otimes 1) \Big( ( \tau \otimes \iota  ) (b_i \otimes 1) \big( ( \tau \otimes \iota  ) \Delta(a) \big) \Big)\bigg)\\
 &= m( \iota  \otimes S )\Big( \sum_{i=1}^{n} ( a_i \otimes 1) ( \tau \otimes \iota  ) \big( ( b_i \otimes 1) \Delta(a) \big) \Big).
\end{align*}
It follows that $\chi(a)b \in A$ for all $a,b \in A$. In addition, for all $c \in A$,
\begin{align*}
    \chi(a)(bc) &= m( \iota \otimes S )\big( ( bc \otimes 1) ( \tau \otimes \iota  ) \Delta(a) \big)\\
    &= m( \iota  \otimes S )\big( ( b \otimes 1 ) ( c \otimes 1) ( \tau \otimes \iota  ) \Delta(a) \big)\\
    &= b \big(m( \iota  \otimes S) ( c \otimes 1 )( \tau \otimes \iota  ) \Delta(a)\big )\\
    &= b \big( \chi(a)c \big).
\end{align*}
This shows that $\chi(a)$ is a right multiplier on $A $ for any $a \in A$. According to the formula $(\ref{3.5})$, the properties of $S$ imply that for any $ a,b \in A $
\begin{align*}
    \chi (a)b \overset{(\ref{3.5})}{=} m(\iota  \otimes S)\Big( ( b \otimes 1) \Delta\big(\tau (a) \big)\Big) = b\varepsilon \big( \tau (a)\big) ,
\end{align*}
this proves that $\chi(a) = \varepsilon \big(\tau( a )\big)$. Therefore, $\chi(A) \subseteq k$ and $\chi$ can be regarded as a map $\chi: A \rightarrow k$. Considering that $\tau$ is an algebra endomorphism and
\begin{align*}
    &\chi (ab) = \varepsilon \big( \tau (ab) \big) = \varepsilon \big( \tau (a) \big) \varepsilon \big( \tau (b) \big) = \chi (a) \chi(b),\\
    &\chi(a + b) = \varepsilon \big( \tau (a + b) \big) = \varepsilon \big( \tau (a) \big) + \varepsilon \big( \tau (b) \big) = \chi(a) + \chi (b),
\end{align*}
which implies that $\chi$ is a character on $A$ and is nondegenerate.

 Now, for $\tau(a) \in A$, equation $(\ref{3.5})$ gives 
\begin{align*}
       ( \iota  \otimes \Delta )\Delta \big( \tau (a) \big) = \tau ( a_{(1)} ) \otimes a_{(2)(1)} \otimes a_{(2)(2)}.
\end{align*}
Then by $(\ref{2.1})$ we have 
\begin{equation}
    S\big( \tau ( a_{(1)} ) \big) a_{(2)(1)} \otimes a_{(2)(2)} = 1 \otimes \tau (a) = \tau( a_{(1)} ) S( a_{(2)(1)} ) \otimes a_{(2)(2)}.  \label{3.7}
\end{equation}
Using the formula $\chi(a) = \varepsilon \big(\tau(a)\big)$, we can express $\tau$ in terms of $\chi$,
\begin{align*}
      \chi (a_{(1)}) a_{(2)}  &= \varepsilon \big( \tau( a_{(1)} ) \big) a_{(2)} = \varepsilon \big( \tau( a_{(1)}) \big) \varepsilon( a_{(2)(1)}) a_{(2)(2)} \\
     &= \varepsilon \big( \tau( a_{(1)}) \big) \varepsilon \big( S(a_{(2)(1)}) \big) a_{(2)(2)}\\
     & = \varepsilon \big(\tau( a_{(1)})S(a_{(2)(1)})\big) a_{(2)(2)}\\
     &\overset{(\ref{3.7})}{=} \varepsilon (1) \tau (a) = \tau (a).
\end{align*}
This establishes equation $(\ref{3.2})$. By substituting $\tau(a) = \chi ( a_{(1)})a_{(2)}$(as established above) into both sides of equations $(\ref{3.5})$ and $(\ref{3.6})$, it ensues that
\begin{align*}
     &\chi(a_{(1)}) a_{(2)(1)} \otimes a_{(2)(2)} = \chi( a_{(1)(1)}) a_{(1)(2)} \otimes a_{(2)} ,\\
    &\chi(a_{(1)}) a_{(2)(1)} \otimes a_{(2)(2)} = A{d_r}( a_{(1)}) \chi(a_{(2)(1)}) \otimes a_{(2)(2)}.
\end{align*}
These two equations are exactly $(\ref{3.3})$.
\smallskip

\textbf{2. Proof of sufficiency.} 

The sufficiency part of the proof is carried out in three steps. 

\textbf{step 1.} Comultiplication. 

Under conditions $(\ref{3.2})$-$(\ref{3.4})$, we proceed to show that $(\ref{3.5})$ and $(\ref{3.6})$ hold. We have 
\begin{align*}
     \Delta \big( \tau ( a ) \big) &\overset{(\ref{3.2})}{=}\Delta \big( \chi ( a_{(1)} ) a_{(2)} \big) = \chi ( a_{(1)}) a_{(2)(1)} \otimes a_{(2)(2)}\\
    & \overset{(\ref{3.3})}{=} \chi(a_{(1)(1)}) a_{(1)(2)} \otimes a_{(2)} \overset{(\ref{3.2})}{=} \tau( a_{(1)}) \otimes a_{(2)},\\
    \Delta \big( \tau ( a ) \big) &\overset{(\ref{3.2})}{=} \Delta \big( \chi ( a_{(1)} ) a_{(2)} \big) = \chi(a_{(1)}) a_{(2)(1)} \otimes a_{(2)(2)}\\
    & \overset{(\ref{3.3})}{=} Ad_r( a_{(1)}) \otimes \chi( a_{(2)(1)}) a_{(2)(2)} \overset{(\ref{3.2})}{=}  Ad_r(a_{(1)}) \otimes \tau(a_{(2)}).
\end{align*}
These equations confirm that $(\ref{3.5})$ and $(\ref{3.6})$ hold.
Together with the necessity part of the proof, we conclude that $\Delta(y)\Delta(a) = \Delta \big(\tau(a)\big)\Delta(y) + \Delta\big(\delta(a)\big)$ is satisfied. The comultiplication $\Delta: A[y; \tau, \delta] \rightarrow M(A[y; \tau, \delta] \otimes A[y; \tau, \delta])$ is defined as above, note that it is not necessarily coassociative.
\smallskip

\textbf{Step 2.} Counit. 

Assume that $A[y; \tau, \delta]$ is equipped with the comultiplication defined in Step 1. Suppose further that the counit $\varepsilon$ can be extended to $A[y; \tau, \delta]$, with $\varepsilon(y) = 0$, and remains an algebra homomorphism. Then the relation $ya = \tau(a)y + \delta(a)$ must be preserved under $\varepsilon$, i.e.,
\begin{align*}
     \varepsilon(y) \varepsilon(a) 
     = \varepsilon \big( \tau (a) \big) \varepsilon(y) + \varepsilon \big(\delta (a) \big)
\end{align*}
for all $a \in A$. Substituting $\varepsilon(y) = 0$, we find that such an extension exists if and only if 
\begin{align*}
    \varepsilon\big(\delta(a)\big)=0.
\end{align*}

We now show that $\varepsilon \big(\delta(a)\big) = 0$ indeed holds for all $a \in A$. Applying the map $m(\iota \otimes \varepsilon)$ to $(\ref{3.4})$, we obtain for any $a \in A$:
\begin{align*}
\delta (a) &= m(\iota  \otimes \varepsilon)\Delta \big( \delta (a) \big) \\
& \overset{(\ref{3.4})}{=} m(\iota  \otimes \varepsilon) \Big((\delta \otimes \iota)\Delta(a) + (r\otimes 1)\big((\iota \otimes \delta)\Delta(a)\big)\Big) \\
&= m( \iota  \otimes \varepsilon)\big( \delta ( a_{(1)} ) \otimes a_{(2)} + ra_{(1)} \otimes \delta (a_{(2)}) \big)\\
&= \delta ( a_{(1)} ) \varepsilon ( a_{(2)} ) + r a_{(1)} \otimes \varepsilon \big( \delta ( a_{(2)} ) \big)\\
&= \delta \big( a_{(1)}\varepsilon ( a_{(2)} ) \big) + r a_{(1)} \otimes \varepsilon \big( \delta ( a_{(2)} ) \big)\\
&= \delta(a) + ra_{(1)} \otimes \varepsilon \big( \delta ( a_{(2)} ) \big).
\end{align*}
Comparing both sides gives $a_{(1)}\varepsilon \big(\delta (a_{(2)})\big) = 0$. Finally,
\begin{align*}
    \varepsilon \big( \delta (a) \big) = \varepsilon \Big( \delta \big( \varepsilon ( a_{(1)} ) a_{(2)} \big) \Big) = \varepsilon ( a_{(1)} ) \varepsilon \big( \delta ( a_{(2)} ) \big) = \varepsilon \Big( a_{(1)} \varepsilon \big( \delta ( a_{(2)} ) \big) \Big) = 0.
\end{align*}
This establishes that $\varepsilon \big(\delta (a)\big) = 0$, confirming the existence of an extension $\varepsilon: A[y; \tau, \delta] \rightarrow k$ that preserves the relation $ya = \tau(a)y + \delta(a)$. For any elements $\sum_{i=1}^{m} a_iy^i, \sum_{j=1}^{n} b_jy^j \in A[y; \tau , \delta ]$, it is straightforward to verify that
\begin{align*}
    (\varepsilon \otimes \iota)\big( \Delta (\sum_{i=1}^{m} a_i y^i)( 1 \otimes \sum_{j=1}^{n} b_j y^j) \big) &= ( \sum_{i=1}^{m} a_i y^i) ( \sum_{j=1}^{n} b_j y^j ),\\
    (\iota \otimes \varepsilon)\big(( \sum_{j=1}^{n} b_j y^j \otimes 1) \Delta( \sum_{i=1}^{m} a_i y^i ) \big) &= (\sum_{j=1}^{n} b_j y^j) (\sum_{i=1}^{m} a_i y^i).
\end{align*}
Thus, we conclude that $\varepsilon: A[y; \tau, \delta] \rightarrow k$ indeed serves as a counit for $A[y; \tau, \delta]$.
\smallskip

\textbf{Step 3.} Antipode.

Suppose $A[y; \tau, \delta]$ is endowed with the comultiplication and counit as defined above. If an antipode $S$ exists on $A[y; \tau, \delta]$, extended from the underlying algebra $A$ by setting $S(y) = -r^{-1}y$, then $S$ satisfies $(\ref{2.1})$ and $(\ref{2.2})$, while also preserving the relation $ya = \tau(a)y + \delta(a)$. This implies that
\begin{eqnarray}
    (m\otimes \iota)(S\otimes \Delta)\Delta(ya) = 1\otimes ya = (m\otimes \iota)(\iota \otimes S\otimes \iota)(\iota \otimes \Delta)\Delta(ya), \label{3.8}\\
    (\iota \otimes m)(\Delta \otimes S)\Delta(ya)= ya\otimes 1 = (\iota\otimes m)(\iota \otimes S\otimes \iota)(\Delta \otimes \iota)\Delta(ya), \label{3.9}
\end{eqnarray}
and 
\begin{align}
    S(a)S(y) = S(y)S\big(\tau (a)\big) + S\big(\delta (a)\big).  \label{3.10}
\end{align}

We now proceed to prove equations $(\ref{3.8})$-$(\ref{3.10})$. For $\delta(a) \in A$ it follows from $(\ref{2.1})$ and $(\ref{3.4})$ that
\begin{align}
     1 \otimes \delta (a)&\overset{(\ref{2.1})}{=}
     (m \otimes \iota )(S \otimes \Delta)\Delta \big(\delta (a )\big)  \nonumber\\
     &\overset{(\ref{3.4})}{=} (m \otimes \iota )(S \otimes \Delta)\Big((\delta \otimes \iota)\Delta(a) + (r\otimes 1)\big((\iota \otimes \delta)\Delta(a)\big)\Big) \nonumber\\
     & = (m \otimes \iota)(S \otimes \Delta )\big(\delta (a_{(1)}) \otimes a_{(2)} + ra_{( 1)} \otimes \delta(a_{(2)})\big)   \nonumber\\
    &  = (m \otimes \iota )\Big(S\big(\delta( a_{(1)})\big)\otimes \Delta ( a_{( 2 )}) + S(ra_{(1)}) \otimes \Delta \big( \delta (a_{( 2 )})\big)\Big)   \nonumber\\
    &  = (m \otimes\iota)\Big(S\big(\delta ({a_{(1)}})\big) \otimes a_{(2)(1)} \otimes a_{(2)(2)}+ S(ra_{(1)}) \otimes \delta (a_{(2)(1)}) \otimes a_{(2)(2)}   \nonumber\\
    &\quad+ S(ra_{(1)}) \otimes ra_{(2)(1)} \otimes \delta( a_{(2)(2)})\Big)   \nonumber\\
    &  = S\big(\delta(a_{(1)})\big)a_{(2)(1)} \otimes a_{(2)(2)} + S( ra_{(1)})\delta(a_{(2)(1)}) \otimes a_{(2)(2)} \nonumber\\
    & \quad+ S(ra_{(1)})ra_{(2)(1)} \otimes \delta(a_{(2)(2)}). 
  \label{3.11}
\end{align}
Applying the formula mentioned above, we deduce that for all $a, b \in A$
\begin{align*}
    &\big(S(b) \otimes 1 \big)(m \otimes \iota )(S \otimes \Delta )\Delta (ya)\\
    &=(m \otimes \iota )(S \otimes \Delta )\big(\Delta (ya)(b \otimes 1)\big)\\
    & = (m \otimes \iota )(S \otimes \Delta )\Big(\Delta \big(\tau (a)y + \delta (a)\big)(b \otimes 1)\Big)\\
    &= (m \otimes \iota )(S \otimes \Delta )\bigg(\Big(\Delta \big(\tau (a)\big)\Delta(y) + \Delta \big(\delta (a)\big)\Big)(b \otimes 1)\bigg)\\
    &\overset{(\ref{3.5})}{=} (m \otimes \iota )(S \otimes \Delta )\Big(\tau (a_{(1)})yb \otimes a_{(2)} + \tau (a_{(1)})rb \otimes {a_{(2)}}y + \Delta \big(\delta (a)\big)(b \otimes 1)\Big)\\
    &= (m \otimes \iota )\Big(S\big(\tau (a_{(1)})yb\big) \otimes a_{(2)(1)} \otimes a_{(2)(2)} + S\big(\tau (a_{(1)})r b\big) \otimes {a_{(2)(1)}}y \otimes a_{(2)(2)} \\
    & \quad +S\big(\tau (a_{(1)})r b\big) \otimes a_{(2)(1)}r \otimes a_{(2)(2)}y \Big)\\
    &\quad+(m \otimes \iota )(S \otimes \Delta )\big(\delta (a_{(1)})b \otimes a_{(2)} + ra_{(1)}b \otimes \delta (a_{(2)})\big)\\
    &\overset{(\ref{3.4})}{=} S\big(\tau (a_{(1)})yb\big){a_{(2)(1)}} \otimes a_{(2)(2)} + S\big(\tau (a_{(1)})r b\big)a_{(2)(1)}y \otimes a_{(2)(2)} \\
    & \quad + S\big(\tau ({a_{(1)}})rb\big)a_{(2)(1)}r \otimes a_{(2)(2)}y\\
    &\quad +(m \otimes \iota)\Big(S\big(\delta (a_{(1)})b\big) \otimes a_{(2)(1)} \otimes a_{(2)(2)} + S(ra_{(1)}b) \otimes \delta (a_{(2)(1)}) \otimes {a_{(2)(2)}} \\
    &\quad + S(ra_{(1)}b) \otimes ra_{(2)(1)} \otimes \delta (a_{(2)(2)})\Big)\\
    &= S(b)S(y)S\big(\tau (a_{(1)})\big)a_{(2)(1)} \otimes a_{(2)(2)} + S(b)S(r) S\big(\tau (a_{(1)})\big){a_{(2)(1)}}y \otimes a_{(2)(2)} \\
    &\quad + S(b)S(r) S\big(\tau (a_{(1)})\big)a_{(2)(1)}r \otimes a_{(2)(2)}y+ S\big(\delta (a_{(1)}b\big)a_{(2)(1)} \otimes a_{(2)(2)} \\
    &\quad + S(ra_{(1)}b)\delta (a_{(2)(1)}) \otimes a_{(2)(2)} + S(ra_{(1)}b)r a_{(2)(1)} \otimes \delta (a_{(2)(1)})\\
    & \overset{(\ref{3.7})}{=} S(b)S(y) \otimes \tau (a) + S(b)S(r)y \otimes \tau(a) + S(b)S(r)r \otimes \tau (a)y \\
    &\quad +(S(b) \otimes 1)\Big(S\big(\delta(a_{(1)})\big)a_{(2)(1)} \otimes a_{(2)(2)} + S(ra_{(1)})\delta (a_{(2)(1)}) \otimes {a_{(2 )(2)}} \\
    &\quad+ S( ra_{(1)})ra_{(2)(1)}\otimes \delta (a_{(2)(2)})\Big)\\
    & \overset{(\ref{3.11})}{=}\big(S(b)S(y) + S(b)S(r)y \big) \otimes \tau (a) + S(b) \otimes \tau (a)y +(S(b) \otimes 1)(1 \otimes \delta(a))\\
    & = S(b) \otimes \tau (a)y +  S(b) \otimes \delta(a)\\
    & = S(b) \otimes \big(\tau (a)y +\delta(a)\big)\\
    & = S(b) \otimes ya\\
    & = \big(S(b) \otimes 1\big)(1 \otimes ya).
\end{align*}
By comparing both sides of the equation above and using the nondegeneracy of the multiplication in $A$, we conclude that $(\ref{3.8})$ holds. Equation $(\ref{3.9})$ can be verified in a similar manner.

To verify condition $(\ref{3.10})$, we make use of the identities $ya = \tau(a)y+\delta(a)$ and $S(y) = -r^{-1}y$  \big(where $r^{-1}y \in M(A[y; \tau, \delta])$\big), which lead to 
\begin{align*}
     -S(a)r^{- 1}y &= -r^{ - 1}yS\big(\tau(a)\big) + S\big(\delta(a)\big) \\
      &=-r^{-1}\tau\Big(S\big(\tau(a)\big)\Big)y - r^{- 1}\delta \Big(S\big(\tau(a)\big)\Big)+S\big(\delta (a)\big)
\end{align*}
for all $a \in A$.
Rearranging terms, we obtain
\begin{align*}
    \bigg(r^{-1}\tau\Big(S\big(\tau(a)\big)\Big)-S(a)r^{- 1}\bigg)y = - r^{- 1}\delta \Big(S\big(\tau(a)\big)\Big)+S\big(\delta (a)\big).
\end{align*}
It follows that condition $(\ref{3.10})$ is equivalent to the following pair of identities holding for all $a \in A$:
\begin{align}
   &S(a)r^{-1}= r^{-1}\tau\Big(S\big(\tau(a)\big)\Big),   \label{3.12}\\
   &\delta\Big(S\big(\tau(a)\big)\Big)=rS\big(\delta(a)\big).      \label{3.13}
\end{align}

It remains to show that $(\ref{3.12})$ and $(\ref{3.13})$ hold for every $a\in A$. Indeed,
\begin{align*}
  \tau\Big( S\big( \tau(a) \big) \Big) &\overset{(\ref{3.2})}{=} \tau\Big(S\big(\chi ({a_{(1)}}){a_{(2)}}\big)\Big) = \tau \big(\chi ({a_{(1)}})S({a_{(2)}})\big) = \chi ({a_{(1)}})\tau \big(S({a_{(2)}})\big)\\
 &\overset{(\ref{3.2})(\ref{3.3})}{=} 
 \chi ({a_{(1)}})A{d_r}\Big(S(a_{(2)(2)})\chi\big(S( a_{(2)( 1 )})\big)\Big)\\
 &= \chi ({a_{(1)}}) \chi\big(S( a_{(2)(1)})\big) A{d_r}\big(S(a_{(2)(2)})\big)\\
 &= \chi \big(a_{(1)}S(a_{(2)(1)})\big)Ad_r\big(S(a_{(2)(2)})\big) \\
 &\overset{(\ref{2.1})}{=} Ad_r\big(S(a)\big),
\end{align*} 
which confirms $(\ref{3.12})$. 

Proceeding to equation $(\ref{3.13})$, we utilize equation $(\ref{3.2})$, from which it follows that
\begin{align*}
    S\big(\tau(a)\big) \overset{(\ref{3.2})}{=} S\big(\chi ({a_{(1)}}){a_{(2)}}\big) = \chi (a_{(1)})S(a_{(2)}).
\end{align*}
Thus, equation $(\ref{3.13})$ can be reformulated equivalently as
\begin{equation}
    \chi (a_{(1)})\delta \big( S(a_{(2)})\big) = rS\big(\delta (a)\big),    \label{3.14}
\end{equation}
it suffices to demonstrate the validity of $(\ref{3.14})$.
Applying the map $m(\iota \otimes S )$  to both sides of the identity
\begin{align*}
   (b \otimes 1) \Delta\big(\delta(a)\big) =b \delta (a_{(1)}) \otimes a_{(2)}+ bra_{(1)} \otimes \delta (a_{(2)}),
\end{align*} 
we derive 
\begin{align*}
    m(\iota \otimes S )\Big((b \otimes 1)\Delta \big(\delta(a)\big)\Big) = b\delta (a_{(1)}) S(a_{(2)}) + bra_{(1)} S(\delta \big({a_{(2)}})\big).
\end{align*}
Using the defining property of the antipode, the left-hand side simplifies as
\begin{align*}
    m(\iota \otimes S )\Big((b \otimes 1)\Delta \big(\delta(a)\big)\Big) = b  \varepsilon\big( \delta (a)\big)= 0,
\end{align*}
which implies
\begin{align*}
       bra_{(1)}S\big(\delta(a_{(2)})\big) = -b\delta(a_{(1)}) S(a_{(2)})  \qquad \text{for any $a,b\in A$}.
\end{align*}  
By the nondegeneracy of the multiplication in $A$, we conclude that the following identity holds for all $a \in A$
\begin{equation}
    a_{(1)}S\big(\delta(a_{(2)})\big)=-r ^{-1}\delta(a_{(1)}) S(a_{(2)}).   \label{3.15}
\end{equation}
Therefore, we have 
\begin{align}
    rS\big(\delta(a)\big) &\overset{(\ref{2.1})}{=} rS(a_{(1)})a_{(2)(1)}S\big(\delta (a_{(2)(2)})\big)  \nonumber\\
    &\overset{(\ref{3.15})}{=}  - rS(a_{(1)})r^{-1} \delta( a_{(2)(1)})S(a_{(2)(2)}) \nonumber\\
   & = -Ad_r \big(S(a_{(1)})\big) \delta(a_{(2)(1)})S(a_{(2)(2)})      \label{3.16}
\end{align}
for all $a \in A$.
Furthermore, applying $\delta$ to the identity $\varepsilon(a) = a_{(1)}S(a_{(2)})$, we find
\begin{align*}
0=\delta(a_{(1)})S(a_{(2)})+\tau(a_{(1)})\delta\big(S(a_{(2)})\big),
\end{align*}
which implies
\begin{equation}
    \tau(a_{(1)})\delta \big(S(a_{(2)}) \big) = -\delta (a_{(1)})S(a_{(2)}).    \label{3.17}
\end{equation}
By combining $(\ref{3.16})$ and $(\ref{3.17})$, we obtain the following
\begin{align*}
    \chi (a_{(1)})\delta \big(S( a_{(2)}) \big) &\overset{(\ref{2.1})}{=} \chi (a_{(1)})\tau \big(S(a_{(2)(1)})a_{(2)(2)(1)}\big)\delta \big(S(a_{(2)(2)(2)}) \big)\\
    &  = \chi (a_{(1)})\tau\big({S(a_{(2)(1)}}) \big)\tau(a_{(2)(2)(1)})\delta \big(S(a_{(2)(2)(2)})\big)\\
    &\overset{(\ref{3.17})}{=}  - \chi (a_{(1)})\tau \big(S(a_{(2)(1)})\big)\delta (a_{(2)(2)(1)})S(a_{(2)(2)(2)})\\
    &  =  - \tau \Big( S\big( \chi (a_{(1)})a_{(2)(1)}\big)\Big)\delta(a_{(2)(2)(1)})S(a_{(2)(2)(2 )})\\
    & \overset{(\ref{3.3})}{=}  - \tau \Big( S\big( \chi (a_{(1)(1)})a_{(1)(2)}\big)\Big)\delta(a_{(2)(1)})S(a_{(2)(2)})\\
    &  =  - \chi (a_{(1)(1)})\tau \big(S(a_{(1)(2)})\big)\delta ( a_{(2)(1)})S(a_{(2)(2)})\\
    &  =  - \chi (a_{(1)(1)})Ad_r\big( S(a_{(1)(2)(2)}) \big)\chi \big( S(a_{(1)(2)(1)})\big)\delta (a_{(2)(1)})S(a_{(2)(2)})\\
    &  =  - \chi \big(a_{(1)(1)}S(a_{(1)(2)(1)})\big) Ad_r\big(S(a_{(1)(2)(2)})\big)\delta(a_{(2)(1)})S(a_{(2)(2)})\\
    & = - Ad_r\big(S(a_{(1)})\big)\delta ( a_{(2)(1)})S(a_{(2)(2)})\\
    &\overset{(\ref{3.16})}{=} r S \big( \delta (a) \big)
\end{align*}
for any $a,b \in A$.
This establishes the validity of $(\ref{3.14})$, and consequently, equation $(\ref{3.13})$ holds for every $a \in A$. It follows that the antipode $S$ can be extended from $A$ to the Ore extension $A[y; \tau, \delta]$.
$\hfill \square$
\\

 \textbf{Remark 3.5.} The proof of Theorem 3.4 adopts a strategy analogous to that of Theorem 3.3 in \cite{JM14} and the theorem in \cite{ZL12}, albeit with notable differences in technical execution. If the underlying algebra of $A$ is unital, then $A$ becomes a Hopf coquasigroup, recovering the framework studied in \cite{JM14}. And if the comultiplication in $A$ is coassociative, then $A$ reduces to a multiplier Hopf algebra, corresponding to the setting treated in \cite{ZL12}. 
\\

 \textbf{Example 3.6.} $(1)$ Let $(A, \Delta, \varepsilon_A, S_A)$ be a multiplier Hopf algebra and $(H, \Delta, \varepsilon_H, S_H)$ a Hopf coquasigroup. Then the tensor product $A \otimes H$ forms a multiplier Hopf coquasigroup under the following structure: 
  \begin{align*}
     & (a \otimes h)(a' \otimes h') = aa' \otimes hh', \\
     & \Delta(a\otimes h) = (a_{(1)} \otimes h_{(1)}) \otimes (a_{(2)} \otimes h_{(2)}) \in M\big((A \otimes H) \otimes (A \otimes H)\big), \\
     & \varepsilon_{A \otimes H}(a \otimes h) = \varepsilon_{A}(a )\varepsilon_{H}(h), \\
     & S_{A \otimes H}(a\otimes h) = S_A(a) \otimes S_H(h),
  \end{align*}
  where $a, a' \in A$ and $h, h' \in H$.
  If $A[y; \tau, \delta]$ is an Ore extension of $A$, then $A[y; \tau, \delta] \otimes H$ is an Ore extension of the multiplier Hopf coquasigroup $A \otimes H$ with 
  \begin{align*}
     \Delta_{A \otimes H} (y \otimes 1_H) = (y \otimes 1_H) \otimes (1_{M(A)} \otimes 1_H) + (r \otimes 1_H) \otimes (y \otimes 1_H).
  \end{align*}

 Similarly, $A \otimes H[y; \tau, \delta]$ is an Ore extension of $A \otimes H$ with 
 $\Delta_{A \otimes H} ( 1_{M(A)} \otimes y) = (1_{M(A)} \otimes y) \otimes (1_{M(A)} \otimes 1_H) + (1_{M(A)} \otimes r) \otimes (1_{M(A)} \otimes y)$.

 $(2)$ Let $A$ be an algebra spanned by elements $\{e_p e_{\alpha} \mid p \in \mathbb{Z},\alpha \in \mathbb{G}\}$ over $\mathbb{C}$, where $\mathbb{G}$ is a (finite or infinite) quasigroup, and the generators $e_pe_{\alpha}$ satisfy the orthogonality relations $(e_pe_{\alpha})(e_q e_{\beta}) = \delta^K_{p,q} \delta^K_{\alpha,\beta} e_pe_{\alpha}$, with $\delta^K$ denoting the Kronecker delta. 
 
 Define a comultiplication $\Delta$ on $A$ by
 \begin{align*}
   \Delta(e_p e_{\alpha}) = \sum\limits_{\beta \in \mathbb{G}} \sum\limits_{k \in \mathbb{Z}} e_k e_{\beta}\otimes e_{p - k}e_{ \beta^{-1}\alpha}.
\end{align*}
 Note that these infinite sums are well-defined in the strict topology on the multiplier algebra, i.e., multiplication with elements of $A$ yields finite sums. One may verify that $(A, \Delta)$ is a regular multiplier Hopf coquasigroup. The counit is given by $\varepsilon(e_p e_{\alpha}) = \delta^K_{p, 0} \delta^K_{\alpha, e}$, and the antipode by $S(e_p e_{\alpha}) = e_{ - p}e_{\alpha^{-1}}$.

Following the construction in Example 3.9 of \cite{ZL12}, choose a nonzero scalar $\lambda \in \mathbb{C} \setminus \{0\}$. 
Define $r \in M(A)$ as $r = \sum\limits_{k \in \mathbb{Z}} \lambda^k e_k$, a character $\chi: A \rightarrow \mathbb{C}$ by $\chi(e_p e_{\alpha}) = \delta^K_{p+2, 0} \delta^K_{\alpha, e}$, an algebra endomorphism $\tau(e_p e_{\alpha}) = e_{p + 2} e_{\alpha}$ and a $\tau$-derivation $\delta: A \rightarrow A$ by $\delta(e_p e_{\alpha}) = 0$. It is straightforward to verify that these definitions satisfy conditions (\ref{3.2})-(\ref{3.4}). Consequently, the Ore extension $A[y; \tau, \delta]$ is established, thereby endowing 
$A$ with the structure of a multiplier Hopf coquasigroup in the sense of Theorem 3.4.
\\

At the end of this section, we introduce the following theorem concerning the $*$-structure of MHC-Ore extensions. In parallel with Definition 2.4 presented in \cite{V94}, we first define the notion of a multiplier Hopf $*$-coquasigroup.
\\

\textbf{Definition 3.7.} If $A$ is a $*$-algebra over $\mathds{C}$, we call $\Delta$ a comultiplication if it is also a $*$-homomorphism. A multiplier Hopf $*$-coquasigroup is a $*$-algebra with a comultiplication, making it into a multiplier Hopf coquasigroup.
\\

Parallel to the proposition in \cite{V94}, we have the following proposition for multiplier Hopf $*$-coquasigroups.

\textbf{Proposition 3.8.} If $A$ is a multiplier Hopf $*$-coquasigroup, then for all $a\in A$, we have 

$(1)$ $\varepsilon(a^*) = \varepsilon(a)^-$ \big(where $\varepsilon(a)^-$ denotes the conjugate of $\varepsilon(a)$\big);
\smallskip 

$(2)$ $S(S(a)^*)^* = a$.

\textbf{Proof.} $(1)$ Define $\tilde{\varepsilon}(a) = \varepsilon(a^*)^-$. Consider the action of $(\tilde{\varepsilon} \otimes \iota)$ on $(1\otimes a)\Delta(b)$
\begin{align*}
    (\tilde{\varepsilon} \otimes \iota)\big((1\otimes a)\Delta(b)\big)
    &= \Big((\varepsilon \otimes \iota)\big(\Delta(b^*)(1\otimes a^*)\big)\Big)^*=(b^*a^*)^*=ab \\ 
    &= (\varepsilon \otimes \iota)\big((1\otimes a)\Delta(b)\big).
\end{align*}
This implies $\tilde{\varepsilon} = \varepsilon$. Substituting the definition of $\varepsilon$, we conclude that $\varepsilon(a) = \varepsilon(a^*)^-$, and hence $\varepsilon(a^*) = \varepsilon(a)^-$.

\smallskip
$(2)$ Suppose $a \otimes b = \sum_{i=1}^{n} \Delta ( a_i) ( 1 \otimes b_i)$. Applying $m(S \otimes \iota)$ to both sides yields
\begin{align*}
    S(a)b &= m (S \otimes \iota) (a \otimes b)
    = \sum_{i=1}^{n} m(S \otimes \iota)\big(\Delta( a_i) (1 \otimes b_i)\big) \\
    &= \sum_{i=1}^{n} \big(m (S \otimes \iota)\Delta (a_i) \big)b_i 
    = \sum_{i=1}^{n} \varepsilon  (a_i)b_i. 
\end{align*}  
Thus, $S(a)b = \sum_{i=1}^{n} \varepsilon (a_i)b_i$. Taking adjoints and applying the flip map gives $b^* \otimes a^*  = \sum_{i=1}^{n}  ( b_i \otimes 1) \Delta^{cop} (a_i^ *)$. Next, applying
$(\iota \otimes S^{-1})$ and multiplying on the right by $1 \otimes c^*$, it follows from Lemma 3.1 in \cite{Y24} that
\begin{align*}
    b^*S^{-1}(a^*)=\sum_{i=1}^{n} b_i^*\varepsilon(a_i^*)c^*.
\end{align*}
Taking adjoints again, we obtain $cS^{-1}(a^*)^*b = cS (a)b$, which implies $S^{-1}(a^*) = S(a)^*$, and therefore $S(S(a)^*)^* = a$.
$\hfill \square$
\\

\textbf{Theorem 3.9.} Let $A$ be a multiplier Hopf $*$-coquasigroup. Then the MHC-Ore extension $A[y; \tau, \delta]$ is also a multiplier Hopf $*$-coquasigroup provided the following conditions are satisfied: 

    $(1)$ $*\tau$ is an involution;
     
    $(2)$ $\delta(\tau + \iota) = 0$ and $\delta * = * \delta$;
    
    $(3)$ $r^* = r$.

\textbf{Proof.} We begin by extending the $*$-operator to $M(A[y; \tau, \delta])$ such that $y^* = y$. For $y \in M(A)[y; \tau, \delta]$ and any $a \in A$, conditions (1) and (2) imply
\begin{align*}
    (ya)^*&=\big(\tau(a)y + \delta(a)\big)^* = y\tau (a)^* + \delta(a)^*\\
    &= \tau \big(\tau (a)^*\big)y + \delta \big( \tau (a)^*\big) + \delta (a)^*
    = a^*y + \delta \big(\tau (a)\big)^* + \delta (a)^*\\
    &= a^*y + \delta \big(\tau(a) + a \big)^*
    = a^*y + \big(\delta(\tau  + \iota)(a) \big)^* \\
    &= a^*y = a^*y^*.
\end{align*}
Applying the result $a^*y = y\tau (a)^* + \delta(a)^*$ derived above, we have for any $a \in A$ and $x \in M(A)$
\begin{align*}
    a^*x^*y &=(xa)^*y = y\tau (xa)^* + \delta (xa)^*\\
    &= y(\tau (x)\tau(a))^* + ( \delta (x)a + \tau(x)\delta(a))^*\\
    &= y\tau (a)^*\tau (x)^* + a^*\delta(x)^* + \delta (a)^*\tau (x) ^*\\
    &= ( y\tau (a)^*+ \delta (a)^*)\tau (x)^*+ a^*\delta (x)^*\\
    &=(a^*y)\tau(x)^*+a^*\delta(x)^*\\
    &= a^*( y\tau (x)^* + \delta (x)^*),
\end{align*}
this implies $x^*y^* = y^*\tau (x)^* + \delta (x)^*$. Hence, the involution $*$ can be extended to a map $*: M(A) \to M(A)$ such that $y^* = y$.
Moreover, for any element $\sum_{i=1}^{n} a_iy^i \in A[y; \tau, \delta]$, its image under the involution satisfies $(\sum_{i=1}^{n} a_iy^i)^* \in A[y; \tau, \delta]$.
It is straightforward to verify that $A[y; \tau, \delta]$ forms a $*$-algebra under the above definition.

Next, we prove that $\Delta$ is a $*$-homomorphism. From condition (3), it follows that
\begin{align*}
    \Delta( y )^* &= ( y \otimes 1 + r \otimes y)^*
     = y^* \otimes 1 + r^* \otimes y^* \\
    & = y \otimes 1 + r \otimes y
     = \Delta(y) \\
     &= \Delta(y^*).
\end{align*}
Since $\Delta$ is an algebra homomorphism, it implies $\Delta(\sum_{i=1}^{n} a_iy^i)^* = \Delta \big((\sum_{i=1}^{n} a_iy^i)^* \big)$.

Furthermore, it is clear that $\varepsilon (y^*) = 0 = \overline {\varepsilon (y)}$, and more generally, $\varepsilon \big((\sum_{i=1}^{n} a_iy^i)^*\big) = \overline {\varepsilon (\sum_{i=1}^{n} a_iy^i)}$, so  $\varepsilon$ is a $*$-homomorphism. 

It remains to prove that $S\big(S(y)^*\big)^* = y$. Using Proposition 3.3 and condition (3), we obtain 
\begin{align*}
    S(y)^* &= (-r^{-1}y)^*= -y^ *(r^{- 1})^*\\
        &= -yr^{-1}= S(y).
\end{align*}
and therefore
\begin{align*}
    S \big(S(y)^*\big) &= S\big(S(y)\big) = S(-r^{-1}y) \\
    &= -S(y)S(r^{-1}) = r^{-1}yr \\
    &= y.
\end{align*}
Thus, $S \big(S(y)^*\big)^* = y^* = y$, which completes the proof.
$\hfill \square$

\section{Isomorphism of Ore extensions}
\def\theequation{\thesection.\arabic{equation}}
\setcounter{equation}{0}

According to Theorem 3.4, an Ore extension of a regular multiplier Hopf coquasigroup is uniquely determined by the parameters $\chi$, $r$, and $\delta$. We therefore denote such an Ore extension as $A[y; \tau, \delta]$ by $A(\chi, r, \delta)$, where $\chi$ is a character on $A$, $r$ is a group-like element in $M(A)$, and $\delta$ is a $\langle \chi, r \rangle$-derivation. The corresponding Ore extension at the level of multiplier coquasigroups is denoted by $M(A)(\chi, r, \delta)$.
\\

Now we can define an isomorphism of MHC-Ore extensions.

\textbf{Definition 4.1.} Two MHC-Ore extensions $A(\chi, r, \delta)$ and $A'(\chi', r', \delta')$ of regular multiplier Hopf coquasigroups $A$ and $A'$ are said to be isomorphic if there is an isomorphism of multiplier Hopf coquasigroups $\phi$: $A(\chi, r, \delta) \rightarrow A'(\chi', r', \delta')$ such that $\phi(A) = A'$. 

Moreover, in the multiplier Hopf $*$-coquasigroup case, $A$ and $A'$ are said to be isomorphic if furthermore $\phi$ is a $*$-homomorphism.
\\

In parallel with \cite[Theorem 4.9]{ZL12}, we have the following theorem for regular multiplier Hopf coquasigroups.

\textbf{Theorem 4.2.} For multiplier Hopf coquasigroups $A$ and $A'$, let their Ore extensions $A(\chi, r, \delta)$ and $A'(\chi', r', \delta')$ be defined with communications: $\Delta (y) = y \otimes 1 + r \otimes y$, $\Delta (y') = y' \otimes 1 + r' \otimes y'$, where $r$ and $r'$ are group-like elements in $A$. Then $A(\chi, r, \delta)$ is isomorphic to $A'(\chi', r', \delta')$ if there is an isomorphism $\phi:A \rightarrow A'$ such that
\begin{align*}
     \phi ( r ) = r', \quad  \tau '\big( {\phi (a)} \big) = \phi \big( {\tau (a)} \big),  \quad \delta'(\phi(a))=\phi(\delta(a))+\phi(\tau(a))d-d\phi(a),
\end{align*}
where $d' \in A'$ such that $\Delta ( d' ) = d' \otimes 1 + r' \otimes d'$.

\textbf{Proof.} Let $\phi(y) = y'+d'$. Then $\phi$ extends from a map $\phi:A \rightarrow A'$ to a homomorphism $\phi: A(\chi, r, \delta) \rightarrow A'(\chi', r', \delta')$. We have    
\begin{align*}
\phi(ya)&=\phi\big(\tau(a)y+\delta(a)\big)=\phi\big(\tau(a)\big)\phi\big(y\big)+\phi\big(\delta(a)\big)\\
&=\phi\big(\tau(a)\big)(y'+d')+\phi\big(\delta(a)\big),\\
\phi(y)\phi(a) &=(y'+d')\phi(a)=\tau'\big(\phi(a)\big)y'+\delta'\big(\phi(a)\big)+d'\phi(a)\\
&=\tau'\big(\phi(a)\big)y'+\phi\big(\delta(a)\big)+\phi\big(\tau(a)\big)d'-d'\phi(a)+d'\phi(a)\\  &=\phi\big(\tau(a)\big)y'+\phi\big(\tau(a)\big)d'+\phi\big(\delta(a)\big)\\
&=\phi\big(\tau(a)\big)(y'+d')+\phi\big(\delta(a)\big),
\end{align*}
for any $a \in A$.
Thus, we get
\begin{align*}
    \phi(ya)=\phi(y)\phi(a).
\end{align*}

Furthermore, we compute:
\begin{align*}
    \Delta\big(\phi(y)\big)&= \Delta (y'+d') = \Delta (y') + \Delta(d') \\
    &= y' \otimes1+r'\otimes y' + d' \otimes 1 + r'\otimes d' \\
    & =(y'+d')\otimes1+r'\otimes(y'+d'),\\
    ( \phi\otimes\phi)\Delta (y) &= (\phi \otimes \phi)(y \otimes 1 + r \otimes y)\\
    &= (y'+d') \otimes 1 + r' \otimes (y' + d'),
\end{align*}
which implies
\begin{align*}
\Delta \big(\phi(y)\big) = (\phi  \otimes \phi)\Delta(y).
\end{align*}

Regarding the counit $\varepsilon$, we obtain
\begin{align*}
    d' = (\varepsilon  \otimes \iota)\big(\Delta (d')\big) = \varepsilon (d') + \varepsilon (r')d' = \varepsilon (d') + d',
\end{align*}
this indicates $\varepsilon (d') = 0$. By Proposition 3.3, we also know $\varepsilon (y) = 0$ and $\varepsilon (y') = 0$.
 Hence, 
\begin{align*}
    \varepsilon \big(\phi (y)\big) = \varepsilon (y' + d') = \varepsilon (y') + \varepsilon (d') = 0,
\end{align*}
which shows that $\varepsilon \big(\phi (y)\big) = \varepsilon (y)$.

Now consider the antipode $S$, we find
\begin{align*}
     0 = \varepsilon (d') = m(S \otimes \iota)\Delta (d') = S(d') + S(r')d' = S(d') + r{'^{ - 1}}d',
\end{align*}
and therefore $S(d') =  - r{'^{ - 1}}d'$. 
Consequently, 
\begin{align*}
    &\phi \big(S(y)\big) = \phi ( - {r^{ - 1}}y) =  - \phi ({r^{ - 1}})\phi (y) =  - r{'^{ - 1}}(y' + d'),\\
    &S\big(\phi (y)\big) = S(y' + d') = S(y') + S(d') =  - {r^{ - 1}}y' - r{'^{ - 1}}d' =- r{'^{ - 1}}(y' + d'),
\end{align*}
it follows that $\phi \big(S(y)\big) = S\big(\phi (y)\big)$.

Hence, we establish that $A(\chi, r, \delta)$ and $A'(\chi', r', \delta')$ are isomorphic within the category of MHC-Ore extensions.
$\hfill \square$
\\

Combining Theorem 3.9 and Theorem 4.2, we obtain the following corollary.

\textbf{Corollary 4.3.} Let $A$ and $A'$ be multiplier Hopf $*$-coquasigroups, with their $*$-structures denoted by $*$, the corresponding Hopf-Ore extensions for $A$ and  
$A'$ are $A(\chi, r, \delta)$ and $A'(\chi', r', \delta')$. If there exists a multiplier Hopf $*$-algebra isomorphism $\phi:A \rightarrow A'$ such that the hypotheses of Theorem 4.2 are satisfied \big(where $(d')^* = d'$\big), and the following conditions hold:

    (1) $*\tau$ is an involution;

    (2) $\delta(\tau + \iota)=0$ and $\delta * = * \delta$; 

    (3) $r^* = r$.
 \\  
Then $A(\chi, r, \delta)$ and $A'(\chi', r', \delta')$ are isomorphic as Hopf-Ore extensions, and they are also isomorphic as multiplier Hopf $*$-coquasigroups.

\textbf{Proof.} By Theorem 3.9 and Theorem 4.2, it suffices to show that $\phi$ preserves the $*$-operation, that is, for any $x \in A (\chi, r, \delta)$, we have $\phi(x^*) = \phi(x)^*$.

By the assumptions of this corollary, $\phi$ is a $ * $-homomorphism on $A$, so for all $a \in A$, $\phi(a^*) = \phi(a)^*$.

Moreover, by Theorem 4.2, we have $\phi(y) = y'+d'$. It follows that
\begin{align*}
    &\phi (y^*) = \phi (y) = y' + d',\\
    &\phi (y)^* = (y' + d')^* = y' + (d')^* = y' + d,
\end{align*}
which implies $\phi (y^*) = \phi (y)^*$, the corollary is proved. 
$\hfill \square$

\section*{Acknowledgements}

The work was partially supported by the Nanjing Agricultural University College Students' Innovative Training Program (No. X2025103070154) China Postdoctoral Science Foundation (No. 2019M651764)
and National Natural Science Foundation of China (No. 11601231).

\end{document}